\documentclass{gtart}

\usepackage{amsmath, amssymb, latexsym, graphicx, epstopdf}
\usepackage{euscript}
\usepackage{subfigure}
\usepackage{epsfig}
\usepackage{psfrag}
\usepackage{amsthm}
\usepackage{setspace}

\theoremstyle{plain}

\newtheorem{thm}{Theorem}
\newtheorem*{thm*}{Theorem}
\newtheorem{lem}[thm]{Lemma}
\newtheorem*{lem*}{Lemma}

\newtheorem*{cor*}{Corollary}

\newtheorem*{cla*}{Claim}
\newtheorem{pro}[thm]{Proposition}
\newtheorem*{pro*}{Proposition}

\newtheorem*{fac*}{Fact}
\newtheorem{que}[thm]{Question}
\newtheorem*{que*}{Question}

\newtheorem*{con*}{Conjecture}

\newtheorem*{rem*}{Remark}

\newtheorem*{rems*}{Remarks}

\newtheorem*{defn*}{Definition}

\newcommand{\C}{\ensuremath{\mathbb{C}}}
\newcommand{\s}{\ensuremath{\mathrm{S}}}
\newcommand{\Z}{\ensuremath{\mathbb{Z}}}
\newcommand{\Q}{\ensuremath{\mathbb{Q}}}
\newcommand{\B}{\ensuremath{\mathrm{B}}}
\newcommand{\pslc}{\mathrm{PSL}_2(\C)}
\newcommand{\slc}{\mathrm{SL}_2(\C)}

\renewcommand{\SS}{\mathcal{S}}
\newcommand{\DD}{\mathcal{D}}
\newcommand{\inv}[1]{{#1}^{-1}}
\newcommand{\brac}[2]{\langle \frac{#1}{#2} \rangle}
\newcommand{\sbrac}[1]{\langle #1 \rangle}
\newcommand{\pie}{\pi_1}
\newcommand{\bound}{\partial}
\newcommand{\tr}{\mathrm{Tr}}
\newcommand{\AND}{\quad \mathrm{and} \quad}

\newcommand{\cc}{\mathcal{C}}

\def\co{\colon\thinspace}

\begin{document}

\title{Not all boundary slopes are strongly detected\\ by the character variety}
\author{Eric Chesebro and Stephan Tillmann}
\begin{abstract}
It has been an open question whether all boundary slopes of hyperbolic knots are strongly detected by the character variety. The main result of this paper produces an infinite family of hyperbolic knots each of which has at least one strict boundary slope that is not strongly detected by the character variety.
\end{abstract}
\maketitle

\section{Introduction}

In 1983 Marc Culler and Peter Shalen introduced a method of constructing essential surfaces in 3--manifolds using the set of representations of the fundamental group into $\slc$. In a nutshell: an ideal point of a curve in the character variety gives a non-trivial action of the fundamental group on a Bass--Serre tree, and this action can be used to construct embedded essential surfaces in a 3--manifold.

Given a compact, orientable, irreducible 3--manifold $M$ with boundary consisting of a single torus, one refers to the boundary slope $r$ of an essential surface associated to an ideal point $\tilde{x}$ of a curve in the character variety as \emph{detected by $\tilde{x}$}. If no closed surface can be associated to $\tilde{x}$, then $r$ is \emph{strongly detected by $\tilde{x}$}, otherwise $r$ is \emph{weakly detected by $\tilde{x}$}. A boundary slope of (an essential surface in) $M$ is termed \emph{detected} if there is some ideal point of a curve in the character variety which detects it. A detected boundary slope is \emph{strongly detected} if it is strongly detected by an ideal point of some curve in the character variety, otherwise it is called \emph{weakly detected}. 

A boundary slope of $M$ is \emph{strict} if it is the boundary slope of an essential surface which is not a fibre or a semi-fibre.  For complements of knots in $S^3$ with the standard framing, the slope $0/1$ is detected by abelian representations and any non--zero boundary slope is strict.

It is shown in \cite{ccgls} that the strongly detected boundary slopes of $M$ are precisely the slopes of the sides of the Newton polygon of the $A$--polynomial. This result implies that (generally speaking) many boundary slopes are detected. The first examples of strict boundary slopes that are not strongly detected were given in 2001 by Schanuel and Zhang \cite{sz}. The manifolds involved are a family of graph manifolds, and the slopes in question are weakly detected \cite{z}. In particular, two interesting questions still remained open.

\begin{que}[Cooper--Long 1996]
  Is every strict boundary slope of the complement of a knot in $S^3$ strongly detected?
\end{que}

\begin{figure}[ht]
\begin{center}
  \includegraphics[width=7cm]{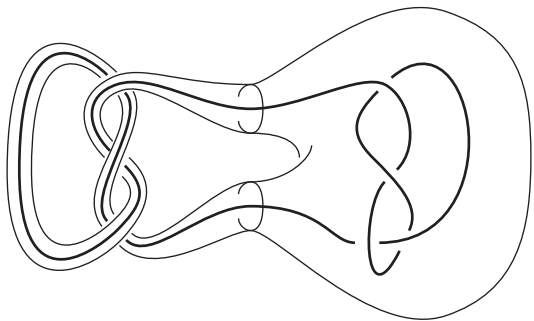}
\end{center}
    \caption{A weakly detected boundary slope}
    \label{fig:connected}
\end{figure}

The following examples of (cylindrical) knot complements containing strict boundary slopes which are weakly detected answer this question negatively.  Consider the connected sum $K$ of two knots in $S^3$ which do not have meridians as boundary slopes: the connected sum of two small knots will suffice. The complement of $K$ contains an essential separating annulus with boundary slope $1/0$. Given the decomposition of the fundamental group of $S^3-K$, it is not difficult to find a curve of characters with ideal points detecting the annulus. However, combining standard arguments involving the limiting representations, any ideal point detecting $1/0$ can be seen to also detect the so--called swallow-follow tori. Thus, the strict boundary slope $1/0$ of $S^3-K$ is weakly detected.

\begin{que}[Schanuel--Zhang 2001]
Let $M$ be an orientable 1--cusped hyperbolic 3--manifold. Is every strict boundary slope of $M$ strongly detected by an ideal point of a curve in its character variety?
\end{que}

\begin{figure}[ht]
\psfrag{a}{{\small $a$}}
\psfrag{b}{{\small $b$}}
\psfrag{c}{{\small $c$}}
\psfrag{d}{{\small $d$}}
\psfrag{e}{{\small $e$}}
\psfrag{f}{{\small $f$}}
\begin{center}
  \includegraphics[width=6cm]{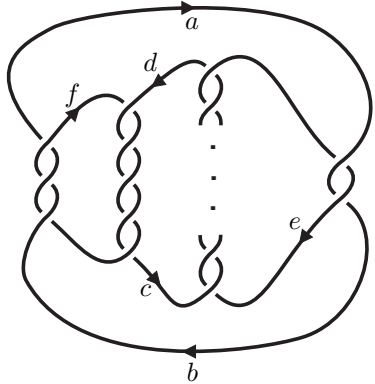}
\end{center}
    \caption{The pretzel knot $K_n=(3,5,2n+1,2)$}
    \label{fig:pretzel}
    \label{fig:Wirt}
\end{figure}

\begin{thm}
The pretzel knot $K_n=(3,5,2n+1,2)$, $n>1$, is hyperbolic and has the strict boundary slope $4(n+4)$ which is not strongly detected by the character variety of $K_n$. 
\end{thm}

The remainder of the introduction discusses the proof of this result; the missing details are given in Sections \ref{sec:noslope} and \ref{sec:strong}. One obvious approach
to proving that a slope is not strongly detected is to calculate the
$A$--polynomial of $K_n$ and to determine the boundary slopes of its
Newton Polygon. This paper takes a different approach --- the key idea is to use 
a relationship between the character varieties of mutants.  It follows from work by Oertel \cite{Oe} that $K_n$ is hyperbolic and that its complement contains two essential 4--punctured spheres (Conway spheres). The knot  $K_n^\tau=(5,3,2n+1,2)$ is obtained from $K_n$ by performing a mutation along the Conway sphere separating the first two tangles from the second two. It is shown in Section \ref{sec:noslope} that an algorithm due to Hatcher and Oertel \cite{HO} implies that $K_n$ has boundary slope $4(n+4)$, but its mutant $K_n^\tau$ does not. 

Work by Cooper and Long \cite{cl96} establishes a relationship between certain strongly
detected boundary slopes of mutants. Let $K$ be a knot in $S^3$, and $K^\tau$ be a Conway mutant of $K$. Denote their character varieties by $X(K)$ and $X(K^\tau)$ respectively. If a curve $C$ in $X(K)$ strongly detects a boundary slope $r$ and there is a character on $C$ whose restriction to the Conway sphere is irreducible, then $r$ is a boundary slope of $K^\tau$ and it is strongly detected by a corresponding curve in $X(K^\tau)$. Details can be found in \cite{cl96} and \cite{tillus_mut}, where it is not explicitly stated that the respective framings are standard. This can be verified with a direct homology argument.

Thus, if the sets of boundary slopes for $K$ and $K^\tau$ are different, then any boundary slope which is strongly detected and not contained in their intersection must be detected by a curve of characters whose restriction to the Conway sphere is reducible. This set of
representations is easier to compute than the whole character variety and, for the knots $K_n$, consists at most of curves detecting the slope $1/0$. The computation can be found in Section \ref{sec:strong} and completes the proof of the theorem.

\begin{rem*}
The above examples also show that not all boundary slopes are strongly detected by the $\pslc$-character variety since each representation into $\pslc$ of the fundamental group of the complement of a knot in $S^3$ lifts to a representation into $\slc$.
\end{rem*}
\begin{que}
Is the boundary slope $4(n+4)$ of $K_n$ weakly detected?
\end{que}
\begin{que}
Is there a small knot complement with a strict boundary slope that is not strongly (and hence also not weakly) detected?
\end{que}

\subsection*{Acknowledgements}

The authors thank Xingru Zhang for sharing with them that the knots $K_3$ and $K_3^\tau$ have different sets of boundary slopes; this was the starting point for this work. Finding the infinite family was made possible through Nathan Dunfield's computer program to compute boundary slopes of Montesinos knots which is freely available at {\tt www.computop.org}\rm.

The second author was supported by a Postdoctoral Fellowship from the Centre de recherches math\'ematiques and the Institut des sciences math\'ematiques in Montr\'eal, and thanks Alan Reid for hosting his visit to the University of Texas at Austin, where this research was conducted.


\section{Mutants with distinct slope sets} \label{sec:noslope}

In this section, the algorithm of \cite{HO} is used to show that the pretzel knot $K_n$ has the boundary slope $4(n+4)$ whereas $K_n^\tau$ does not. For the remainder of this paper, the knots will be denoted in Montesinos' notation by $K_n=K(\frac{1}{3}, \ \frac{1}{5}, \ \frac{1}{2n+1}, \ \frac{1}{2})$ and $K^\tau_n = K(\frac{1}{5}, \ \frac{1}{3}, \ \frac{1}{2n+1}, \ \frac{1}{2})$. This is the notation used in \cite{HO}. The knot $K_n$ is shown in Figure \ref{fig:pretzel}. The knot $K^\tau_n$ is a mutant of $K_n$ by an involution $\tau \co \s_4 \longrightarrow \s_4$, where $\s_4$ is the $4$-punctured sphere which separates the $\frac{1}{3}$ and $\frac{1}{5}$ tangles from the other two tangles.

\subsection{Setup from \cite{HO}} \label{sec:setup}

The constructions from \cite{HO} which apply to our setting are reviewed in this section. In particular, it is assumed throughout that $K$ denotes a four--tangle pretzel knot.

The 3--sphere $\s^3$ can be decomposed into the union of four $3$-balls, $\left\{ \B_i \right\}_{i=0}^3$, where $\bigcap_0^3 \B_i \cong \s^1$ (the {\it axis} for $K$), and each $\B_i$ contains exactly one of the four given tangles.  If $F$ is a properly embedded surface in $M:=\s^3-N(K)$, we may assume that it is transverse to $\bound \B_i$ for every $i$. Hence $F_i := F \cap \B_i$ is properly embedded in $\B_i$ for all $i$ and $P_i:= \bound B_i \cap K$ consists of four points.  In particular, $\bound F_i$ is a curve system on the $4$-punctured sphere $\s_i:=\bound \B_i -N(P_i)$.  We will denote the $i^\mathrm{th}$ tangle $K \cap \B_i$ by $T_i$.

\begin{figure}[t]
\psfrag{a}{{\small $a$}}
\psfrag{b}{{\small $b$}}
\psfrag{c}{{\small $c$}}
\begin{center}
\subfigure[The traintrack with projective weights $a$, $b$, and $c$]{
\includegraphics[width=4.5cm]{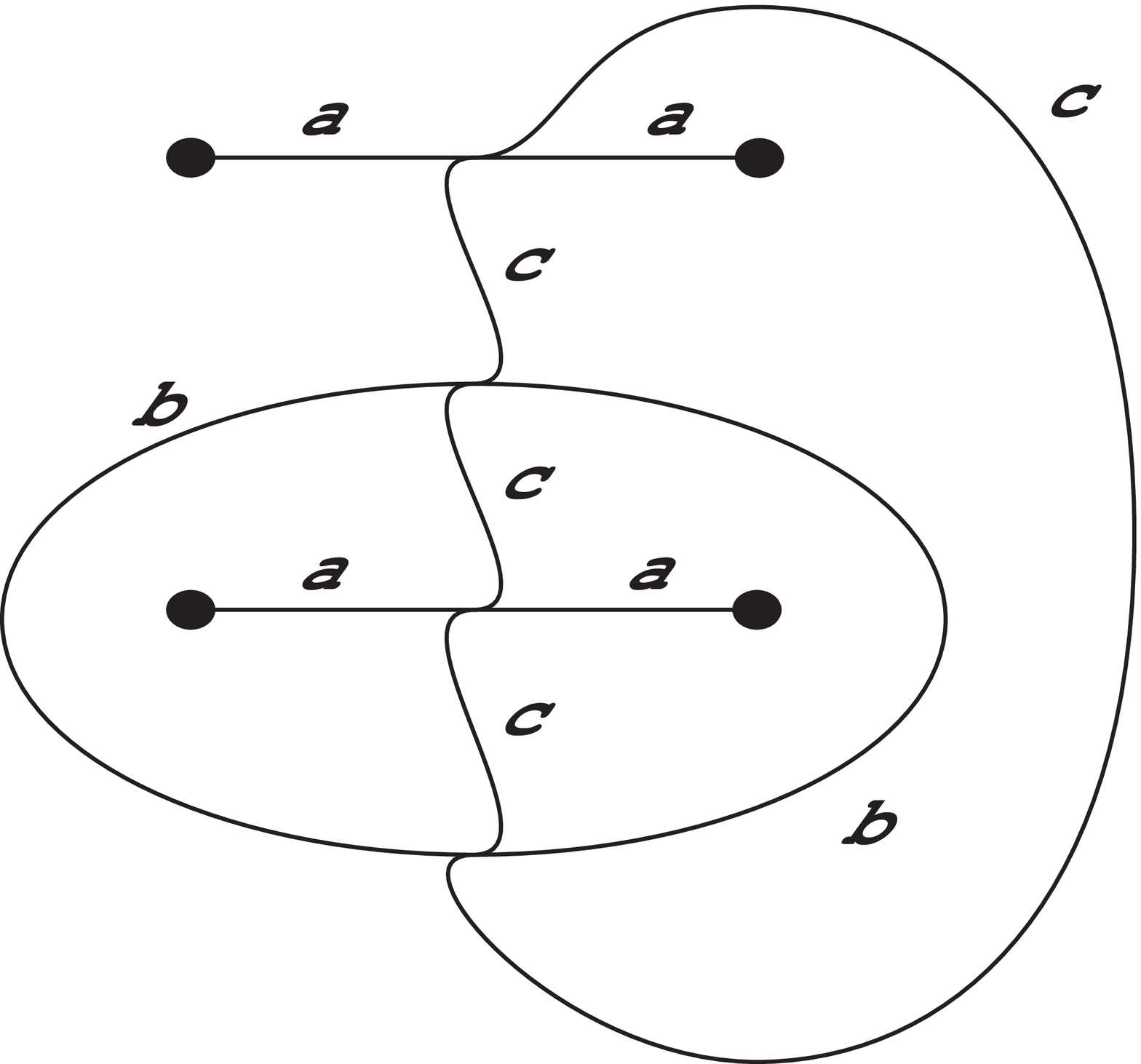}
\label{fig:traintrack}}
\qquad\qquad\qquad
\subfigure[The $\infty$-tangle]{
\includegraphics[width=3cm]{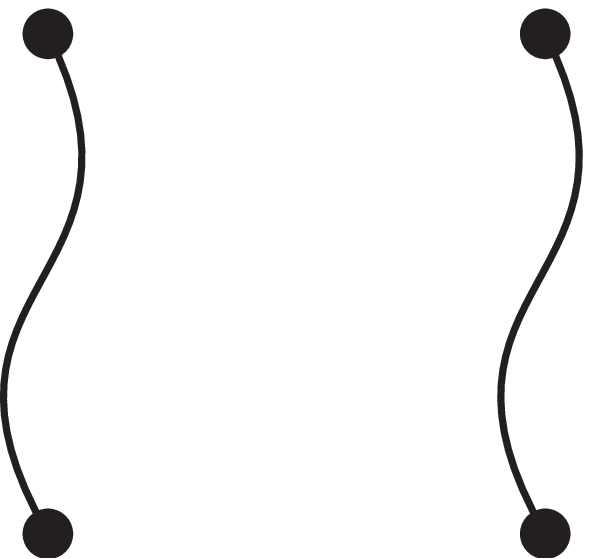}
\label{fig:infty}}
\end{center}
\caption{Curve systems and tangles}  
\end{figure}

Curve systems on the $4$-punctured sphere are either carried by the train track shown in Figure \ref{fig:traintrack} or by its mirror image.  Usually we will consider curve systems up to projective class, so rational projective coordinates $[a,b,c] \in \Q P^2$ represent a projective curve system according to the figure.  We will also use $uv$-coordinates for projective curve systems, where $[a,b,c]$ corresponds to $\left( \frac{b}{a+b}, \, \frac{c}{a+b} \right)$ in $uv$-coordinates.  The slope of the curve system is the $v$-coordinate.

A {\it $p/q$-tangle}, denoted by $\brac{p}{q}$, is a projective curve system $[1,q-1,p]$, or equivalently $( (q-1)/q, \, p/q)$ in $uv$-coordinates.  A {\it $p/q$-circle}, denoted by $\brac{p}{q}^\circ$, is a projective curve system $[0,q,p]$, or equivalently $( 1, \, p/q)$.  The {\it $\infty$-tangle}, denoted by $\langle \infty \rangle$, is the projective class of the pair of vertical arcs shown in Figure \ref{fig:infty}, and will be represented by $(-1,0)$ in $uv$-coordinates.

\begin{figure}[ht] 
\begin{center}
\includegraphics[width=4in]{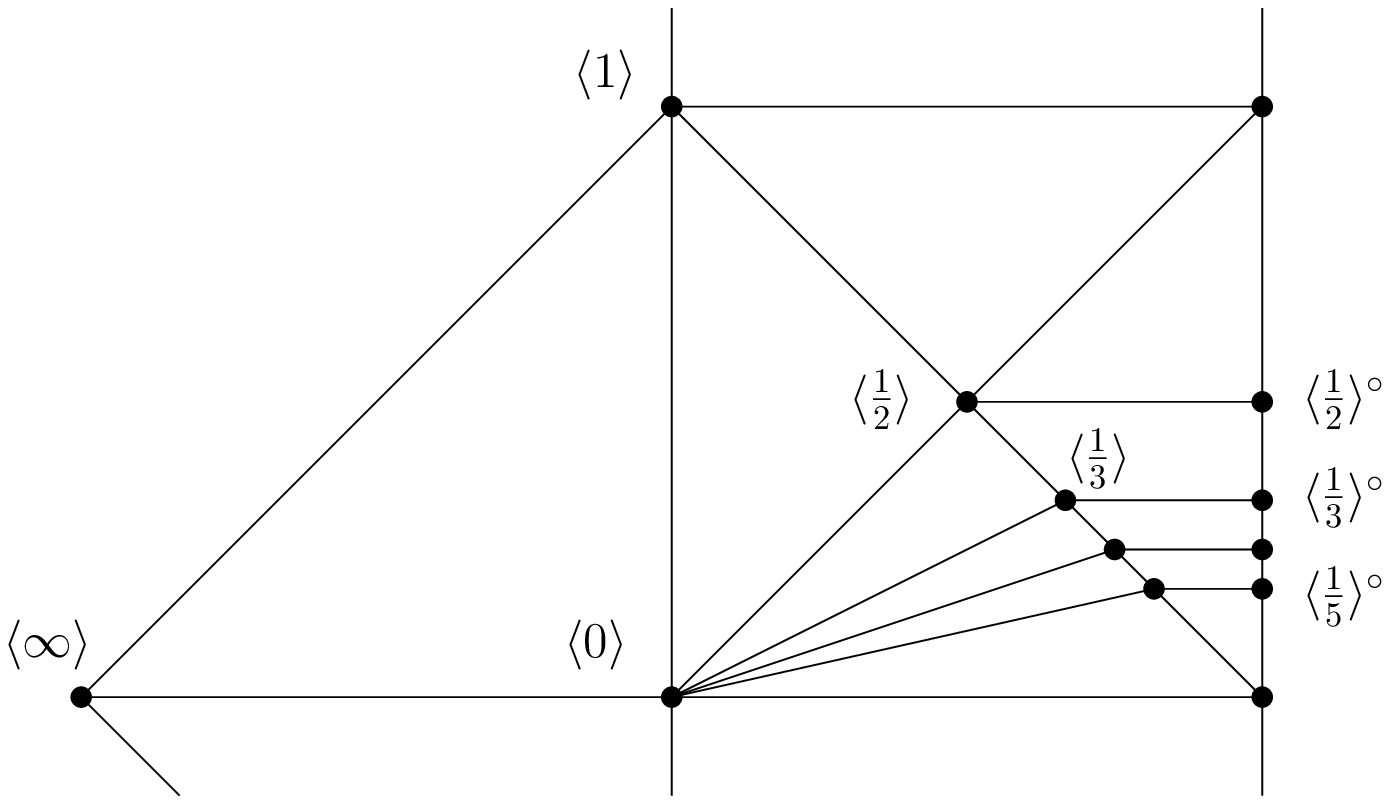}
\end{center}
\caption{A piece of $\DD$} \label{fig:DD}
\end{figure}

We now use the above definitions to define a graph $\DD$ in the $uv$-plane.  The vertices of $\DD$ are the $uv$-coordinates of the $p/q$-tangles and $p/q$-circles for every $p/q \in \Q$ together with the point $\langle \infty \rangle = (-1,0)$.  There are four types of edges in $\DD$, {\it non-horizontal edges, horizontal edges, vertical edges} and {\it infinity edges}.  Two vertices $\brac{p}{q}$ and $\brac{r}{s}$ are connected by a non-horzontal edge if $|ps-qr|=1$ or equivalently if $\brac{r}{s}$ can be obtained from $\brac{p}{q}$ by surgery on an arc.  The horizontal edges connect the vertices $\brac{p}{q}^\circ$ to $\brac{p}{q}$.  The vertical edges connect $\langle n \rangle$ to $\langle n+1 \rangle$ for every $n \in \Z$.  Finally, the infinity edges connect the integer vertices $\langle n \rangle$ to $\langle \infty \rangle$.  If $v$ and $w$ are vertices of $\DD$ that are connected by an edge, we will denote the edge by $[v,w]$.  The subgraph $\SS \subset \DD$ is defined as the portion of $\DD$ with $u$-coordinate in the interval $[0,1]$.  Figure \ref{fig:DD} shows part of the the graph $\DD$.

Given an edge $[v,w]$ in $\DD$, we subdivide it as follows.  For each $m \in \Z^+$ and $k = 1, \ldots , m-1$, let $\frac{k}{m} \cdot v + \frac{m-k}{m} \cdot w$ denote the point on $[v,w]$ corresponding to the curve system consisting of $k$ parallel copies of the pair of arcs representing $v$, together with $m-k$ copies of the pair of arcs representing $w$.

Let $\bound \B_i \times [0,1]$ be a collar on $\bound \B_i$ inside $\B_i$ with $\bound B_i=\bound \B_i \times\{1\}$.  Let $1/q_i$ be the rational number corresponding to the tangle $T_i$.  Then $T_i$ is isotopic (rel $\bound \B_i$) to the two component representative for $\brac{1}{q_i}$ at the level $S_i \times \{0\}$ together with the four arcs $P_i \times [0,1]$.  Our goal is to associate a surface in $\B_i- N(T_i)$ to certain paths in $\DD$.

An {\it edgepath system}, $\gamma=(\gamma_0, \ldots, \gamma_3),$ is an $4 $-tuple of edgepaths in $\DD$ which start and end at rational points of $\DD$.  An {\it admissible} path system is a path system with the following four properties:
\begin{itemize}
\item[(E1)] The starting point of $\gamma_i$ lies on the horizontal edge connecting $\brac{1}{q_i}^\circ$ to $\brac{1}{q_i}$ and if the starting point is not the vertex $\brac{1}{q_i}$, then the path $\gamma_i$ is constant.
\item[(E2)] $\gamma_i$ is {\it minimal}.  That is, it never stops and retraces itself and it never travels along two sides of a triangle in $\DD$ in succession.
\item[(E3)] The ending points of the $\gamma_i$'s all have the same $u$-coordinates and their $v$-coordinates sum to zero.
\item[(E4)] Each $\gamma_i$ proceeds monotonically from right to left, in the sense that traversing vertical edges is permitted.
\end{itemize}  
Admissible edgepath systems are divided into the following three types:
\begin{itemize}
\item A {\it type I system} is an edgepath system, $\gamma$, where each $\gamma_i$ stays in $\SS$ and has no vertical edges.
\item A {\it type II system} is the same as a type I system except that at least one $\gamma_i$ has a vertical edge.
\item A {\it type III system} is a system where the $\gamma_i$'s end to the left of $\SS$.
\end{itemize}

Hatcher and Oertel show how to associate the so--called \emph{candidate surfaces} to each of these types of edgepaths. It is also shown in \cite{HO} that every essential surface in $M$ with non-empty boundary of finite slope is isotopic to such a surface. For the purpose of this paper, it suffices to describe the construction of the candidate surfaces for type II and III systems with no constant paths and with the endpoints of the $\gamma_i$'s on vertices of $\DD$. Let $\gamma$ be such a system. We will write $\gamma_i$ as $[ v_n, \ldots , v_0 ]$ where the $v_j$'s are the vertices of $\gamma_i$ and the vertex $v_j$ is followed by the vertex $v_{j+1}$ as we move along the path.  A complete list of candidate surfaces for $\gamma$ is described in the following paragraphs.

\begin{figure}[ht]
\begin{center}
\includegraphics[width=1.75in]{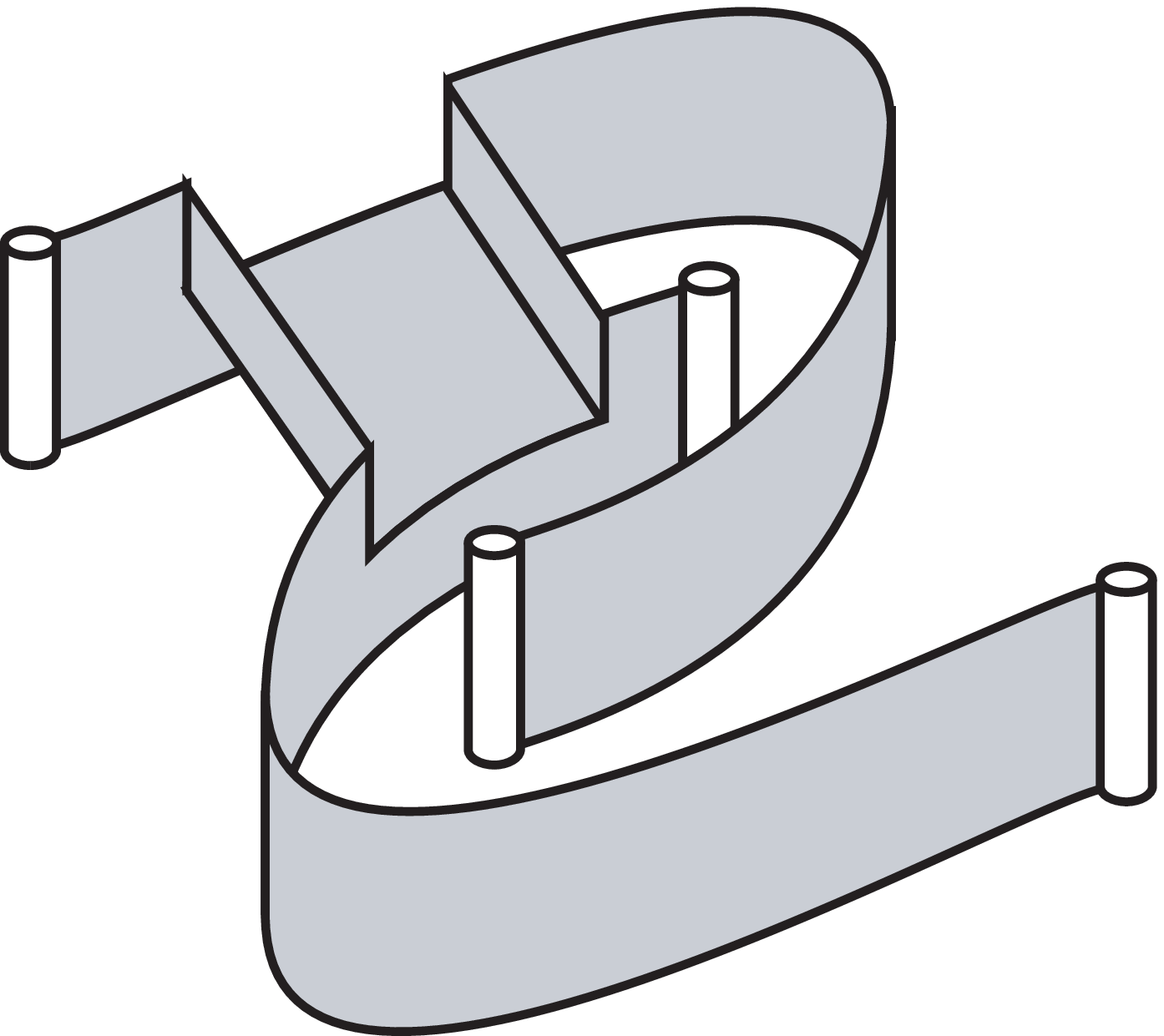}
\end{center}
\caption{A surface corresponding to the edge $\left[\langle1\rangle, \brac{1}{2}\right]$}  \label{fig:edge}
\end{figure}

For any edge $[v_{j+1},v_j]$ of $\gamma_i$ and $a < b,$ we can build a $1$-sheeted surface in $S_i \times [a,b]$ as follows.  Let $\alpha$ and $\beta$ be the curve systems with two arc components in $S_i$ that represent $v_j$ and $v_{j+1}$ respectively.  The surface is 
$$\alpha \times \Big[ a, \, (a+b)/2 \Big) \ \cup \ \beta \times \Big( (a+b)/2, \, b \Big] \ \cup \ D,$$
where $D$ is a regular neighborhood (saddle) of a surgery arc in $S_i \times \left\{(a+b)/2 \right\}$ given by the existence of the the edge $[v_{j+1},v_j]$.  Up to level preserving isotopy there are two choices for each saddle.   One of the two possible surfaces for $\left[\langle1\rangle, \brac{1}{2}\right]$ is shown in Figure \ref{fig:edge}.  To each edgepath $\gamma_i$ we can now associate a surface $\Gamma_i$ which is properly embedded in $\B_i-N(T_i)$ by simply stacking the surfaces associated to the edges of $\gamma_i$ in the collar $S_i \times [0,1]$. Then 
$\Gamma_i \cap (S_i \times \{1\})$
is the two component representative of the vertex $v_n$ and 
$\Gamma_i \cap (S_i \times \{0\})$ 
is the two component representative of the vertex $v_0$ and so lies on the tangle $T_i$.  The condition (E3) guarantees that the surfaces $\{ \Gamma_i \}$ fit together to give a surface in $M$.  

We have described the construction of a $1$-sheeted candidate surface.  If $m$ is a positive integer, we can also construct an $m$-sheeted surface for $\gamma$ in a similar way.  Again, for every edge $[v_j, v_{j+1}]$ in $\gamma_i$, we construct a surface in $S_i \times [a,b]$.  This time we further subdivide the edge into $m$ edges separated by the points $\frac{k}{m} \cdot v_j + \frac{m-k}{m} \cdot v_{j+1}$ where $k = 1, \ldots , m-1$.  Starting with the product surface $\alpha \times \big[a,(a+b)/2m\big)$, we add a saddle at the level $(a+b)/2m$ to pass to the curve system $\frac{1}{m} \cdot v_j + \frac{m-1}{m} \cdot v_{j+1}$.  An $m$-sheeted surface for $[v_j, v_{j+1}]$ is completed by adding saddles to pass between the product surfaces that correspond to each point $\frac{k}{m} \cdot v_j + \frac{m-k}{m} \cdot v_{j+1}$.  Just as before, we can build a surface for $\gamma$ by stacking the surfaces for each edge $\gamma_i$ and finally gluing everything together to give a surface in $M$.   The meridian of the knot will intersect an $m$-sheeted surface $m$ times.  It is worth noticing that, because every time we pass a point $\frac{k}{m} \cdot v_j + \frac{m-k}{m} \cdot v_{j+1}$ we have a choice between two possible saddles, there can be many possible $m$-sheeted candidate surfaces associated to a given path system.

The question whether a candidate surface is essential or not can be studied via the path systems which carry them. This motivates the following terminology. An admissible path system $\gamma$ is {\it incompressible} if every candidate surface associated to $\gamma$ is incompressible, and it is {\it compressible} if every associated candidate surface is compressible. If there are both compressible and incompressible candidate surfaces associated to $\gamma$, it is said to be {\it indeterminant}.

For a path $\gamma_i$ in $\gamma$, let $e_+$ be the number of edges of $\gamma_i$ which increase slope, and $e_-$ be the number of slope decreasing edges. The numbers $e_+$ and $e_-$ are independent of the behavior of $\gamma_i$ to the left of $u=0$.  The {\it twist number} of $\gamma_i$ is $t(\gamma_i)=2(e_- - e_+)$, and the twist number of $\gamma$ is $t(\gamma) = \sum t(\gamma_i)$.  If $F$ is a surface carried by $\gamma$, the twist number $t(F)$ of $F$ is defined to be $t(\gamma)$.  Hatcher and Oertel show that if $S$ is the Seifert surface for a knot and $s$ a path system which carries $S$, then the boundary slope for any surface carried by $\gamma$ is $t(\gamma)-t(s)$.

In order to decide whether an admissible path system is compressible, incompressible, or indeterminant, Hatcher and Oertel define the notions of $r$-values and completely reversible paths.  The {\it $r$-value} of a leftward directed non-horizontal edge is the denominator of the $v$ coordinate of the point where the extension of the edge meets the right hand edge of $\DD$.  The $r$-value is taken to be positive if the edge travels upwards and negative if it travels downwards.  If $\gamma_i$ is a path in $\DD$ then its {\it final $r$-value} is the $r$-value for the last edge in the path.  If $\gamma$ is an admissible path system, then its {\it cycle of final $r$-values} is the $4$-tuple of final $r$-values for the four paths $\{\gamma_i\}$.  The cycle of final $r$-values for $\gamma$ is defined up to cyclic permutation. A path $\gamma_i$ is {\it completely reversible} if each pair of successive edges in $\gamma_i$ lies in triangles that share a common edge.

\subsection{Boundary slopes for $\mathbf{K_n}$ and $\mathbf{K_n^\tau}$}

The goal of this section is to apply the machinery of the previous section to prove that $K_n$ has the boundary slope $4(n+4)$ whereas $K_n^\tau$ does not.  

Consider the admissible path system $s=(s_i)$ given by
\begin{eqnarray*}
s_0 & = & \left[\sbrac{\infty}, \, \sbrac{1}, \, \sbrac{1/2}, \, \sbrac{1/3}\right]\\
s_1 & = & \left[\sbrac{\infty}, \, \sbrac{1}, \, \sbrac{1/2}, \, \sbrac{1/3}, \, \sbrac{1/4}, \, \sbrac{1/5}\right]\\
s_2 & = & \left[\sbrac{\infty}, \, \sbrac{1}, \, \sbrac{1/2},  \, \ldots ,  \, \sbrac{1/(2n+1)}\right]\\
s_3 & = & \left[\sbrac{\infty}, \, \sbrac{1}, \, \sbrac{1/2}\right].
\end{eqnarray*}
Also let $s^\tau=(s_1, \, s_0, \, s_2, \, s_3)$.

\begin{lem} \label{lemF}
A Seifert surface $\Sigma$ for $K_n$ is carried by the path system $s$.  Similarly, a Seifert surface $\Sigma^\tau$ for $K_n^\tau$ is carried by $s^\tau$.  Furthermore, the twist numbers $t(\Sigma)$ and $t(\Sigma^\tau)$ are both $-(14+4n)$.
\end{lem}

\begin{proof} 
While building a $1$-sheeted candidate surface for $s,$ according to the construction from Section \ref{sec:setup}, we may choose saddles so that the resulting surface is is the same as the Seifert surface obtained by applying Seifert's algorithm to the projection of $K_n$ shown in Figure \ref{fig:Wirt}.  We calculate the twist number using the formula $t(\gamma) = \sum t(\gamma_i)$.  We have
\begin{eqnarray*}
t(s_0) & = & 2(0-2) \ = \ -4, \\
t(s_1) & = & 2(0-4) \ = \ -8, \\
t(s_2) & = & 2(0-2n) \ = \ -4n, \\
t(s_3) & = & 2(0-1) \ = \ -2, 
\end{eqnarray*}
and hence
$$t(s) \ = \ \sum t(s_i) \ = \ -4-8-4n-2 \ = \ -(14+4n).$$
The same procedure works for $K^\tau_n$.
\end{proof}

\begin{lem} \label{lem:Ktau_doesn't}
$4(n+4)$ is not a boundary slope of $K^\tau_n$.
\end{lem}

\begin{proof} 
It suffices to show that every candidate surface with twist number $2$ is compressible; this will be proved using three claims which will be established below.

{\bf Claim 1} If $\delta$ is a type III path system with $t(\delta)=2$ then 
\begin{eqnarray*}
\delta_0 & = & \left[ \sbrac{\infty}, \, \sbrac{0}, \, \sbrac{1/5} \right] \\
\delta_1 & = & \left[ \sbrac{\infty}, \, \sbrac{1}, \, \sbrac{1/2}, \, \sbrac{1/3} \right] \\
\delta_2 & = & \left[ \sbrac{\infty}, \, \sbrac{0}, \, \sbrac{1/(2n+1)} \right] \\
\delta_3 & = & \left[ \sbrac{\infty}, \, \sbrac{0}, \, \sbrac{1/2} \right] 
\end{eqnarray*}

Assume $\delta$ is as claimed.  Proposition 2.5 of \cite{HO} implies that since the sum of integer vertices of $\delta$ is $1$, if two of the paths $\delta_i$ are completely reversible then $\delta$ is a compressible path system. It can be verified that $\delta_0$ is completely reversible by noticing that the triangle $[\sbrac{1/5}, \sbrac{0}, \sbrac{1}]$ shares an edge with the triangle $[\sbrac{0}, \sbrac{\infty}, \sbrac{1}]$. The same argument applies to $\delta_2$ and $\delta_3$, hence $\delta$ is a compressible path system.

{\bf Claim 2}  If $\delta$ is a type II path system with $t(\delta)=2$ then (up to adding or removing vertical edges on the individual paths)
\begin{eqnarray*}
\delta_0 & = & \left[ \sbrac{0}, \, \sbrac{1/5} \right] \\
\delta_1 & = & \left[ \sbrac{0}, \, \sbrac{1}, \, \sbrac{1/2}, \, \sbrac{1/3} \right] \\
\delta_2 & = & \left[  \sbrac{0}, \, \sbrac{1/(2n+1)} \right] \\
\delta_3 & = & \left[ \sbrac{0}, \, \sbrac{1}, \,  \sbrac{1/2} \right] 
\end{eqnarray*}

If $\delta$ is as in Claim 2, then the cycle of final $r$-values for $\delta$ is $(-4, \, 1, \, -2n, \, 1)$.  Proposition 2.9 of \cite{HO} shows that there is no incompressible surface associated to $\delta$.

{\bf Claim 3}  There are no type I path systems with twist number $2$.

Assuming the claims, the proof of the lemma is completed. 
\end{proof}

A {\it basic path} is a path starting at a vertex $\langle \frac{p}{q} \rangle$ that proceeds monotonically to the left (without vertical edges) ending at the left edge of $\SS$.  A {\it basic path system} is a path system made up of basic paths.  Note that a basic path $\delta_i$ from $\langle \frac{1}{r} \rangle$  is either $\left[ \sbrac{1}, \, \sbrac{1/2}, \, \sbrac{1/3}, \, \ldots, \, \sbrac{1/r} \right]$ or $\left[ \sbrac{0}, \, \sbrac{1/r} \right]$.  Hence any such path ends at either $0$ or $1$.  If it ends at $0$ then $t(\delta_i)=2(1-0)=2$.  If it ends at $1$ then $t(\delta_i)=2\big(0-(r-1)\big)=2-2r$.  If $\delta$ is a basic path system for $K_n^\tau$ that satisfies (E1) and (E2), then $\delta_0$ starts at $\langle \frac{1}{5} \rangle$, $\delta_1$ starts at $\langle \frac{1}{3} \rangle$, $\delta_2$ starts at $\langle \frac{1}{2n+1} \rangle$, and $\delta_3$ starts at $\langle \frac{1}{2} \rangle$.  So we have 
\begin{eqnarray*}
t(\delta_0) & = & \left\{ \, \begin{array}{rcl} 2 & & \mathrm{if} \ \delta_0 \ \mathrm{ends \ at} \ 0 \\
-8 & &  \mathrm{if} \ \delta_0 \ \mathrm{ends \ at} \ 1 \end{array} \right. \\
t(\delta_1) & = & \left\{ \, \begin{array}{rcl} 2 & & \mathrm{if} \ \delta_1 \ \mathrm{ends \ at} \ 0 \\
-4 & &  \mathrm{if} \ \delta_1 \ \mathrm{ends \ at} \ 1 \end{array} \right. \\
t(\delta_2) & = & \left\{ \, \begin{array}{rcl} 2 & & \mathrm{if} \ \delta_2 \ \mathrm{ends \ at} \ 0 \\
-4n & &  \mathrm{if} \ \delta_2 \ \mathrm{ends \ at} \ 1 \end{array} \right. \\
t(\delta_3) & = & \left\{ \, \begin{array}{rcl} 2 & & \mathrm{if} \ \delta_3 \ \mathrm{ends \ at} \ 0 \\
-2 & &  \mathrm{if} \ \delta_3 \ \mathrm{ends \ at} \ 1 \end{array} \right.
\end{eqnarray*}

\begin{proof}[Proof of Claim 1]
Let $\delta$ be a type III path system with $t(\delta)=2$.  Notice that (E2) implies that $\delta$ has no vertical edges.  Let $\delta'$ be the basic path system $\delta \cap \SS$.  Since extending paths to $\langle \infty \rangle$ doesn't affect twist number, we have $2=t(\delta)=t(\delta')$.  Since there are two choices for each path and four paths, we have a total of $2^4$ possible path systems $\delta'$.  The corresponding twist numbers are in the set
\begin{multline} \nonumber
\{-12, \, -8, \, -6, \, -2, \, 2, \, 4, \, 8, \, 6-4n, \, 2-4n, \, -4n, \\
-4-4n, \, -8-4n, \, -10-4n,  \, -14-4n.\}
\end{multline}
Since $n$ is an integer, $2 \notin \{ -4n, \, -4-4n, \, -8-4n \}$  and $n \geq 2$ so $2 \notin \{ 6-4n, \, 2-4n, \, -10-4n, \, 14-4n \}$.  Then we must have $t(\delta')=2$ for every $n$.  There is only one such path system.  It extends to the system $\delta$ given in the claim.
\end{proof}

\begin{proof}[Proof of Claim 2]
Assume that $\delta$ is a type II path system with $t(\delta)=2$. Then $t(\delta)$ is determined by the basic path system $\delta'$ obtained by deleting all vertical edges of $\delta$. There are infinitely many extensions of $\delta'$, which satisfy (E3), formed by adding vertical edges to the end of the individual paths.  However, even if the extensions do not satisfy the minimality condition (E2), any two such extensions will always have the same twist number.  Thus, we may chose to work with the extension $\tilde{\delta}$ where all paths end at $\langle 0 \rangle$.  Then 
$$t(\tilde{\delta}_i) \ = \ \left\{ \, \begin{array}{ccl} 2 & & \mathrm{if} \ \delta'_i \ \mathrm{ends \ at} \ 0 \\
t(\delta'_i)+2 & &  \mathrm{if} \ \delta'_i \ \mathrm{ends \ at} \ 1 \end{array} \right. .$$
Again there are $2^4$ possibilities.  The corresponding twist numbers are in the set
\begin{multline} \nonumber
\{-6, \, -4, \, -2, \, 0, \, 2, \, 4, \, 6, \, 8, \, 8-4n, \, 6-4n, \,4-4n, \\
 2-4n, \, -4n,  \, -2-4n, \, -4-4n, \, -6-4n.\}
\end{multline}

As before, we know that $2 \notin \{-4-4n, \, -4n, \, 4-4n, \, 8-4n\}$ since $n$ is an integer.  Also, $n \geq 2$ implies $2 \notin \{-6-4n, \, -2-4n, \, 2-4n, \, 6-4n\}$.   Again, the only remaining possibility satisfies the claim.
\end{proof}

\begin{proof}[Proof of Claim 3]
Assume that $\delta$ is a type I path system. Let $( t_0,  \ldots, t_3)$ be the $4$-tuple of endpoints of the paths $\delta_i$.  Since $\delta$ is an admissible path system we know that the sum of the vertical coordinates of the $t_i$'s is zero, but every point of every $\delta_i$ has vertical coordinate greater than or equal to zero.  Therefore $t_i=0$ for every $i$.  We now know that 
\begin{eqnarray*}
\delta_0 & = & \left[ \sbrac{0}, \, \sbrac{1/3} \right] \\
\delta_1 & = & \left[ \sbrac{0},  \, \sbrac{1/5} \right] \\
\delta_2 & = & \left[  \sbrac{0}, \, \sbrac{1/(2n+1)} \right] \\
\delta_3 & = & \left[ \sbrac{0}, \,  \sbrac{1/2} \right] .
\end{eqnarray*}
However, for this path system, we have
$$t(\delta) \ = \ 2+2+2+2 \ = \ 8.$$
We have now established claim 3.
\end{proof}

In order to show that $K_n$ has boundary slope $4(n+4)$, consider the path system $\gamma=(\gamma_i)$ in $\SS$ given by
\begin{eqnarray*}
\gamma_0 & = & \left[\sbrac{1}, \, \sbrac{1/2}, \, \sbrac{1/3}\right]\\
\gamma_1 & = & \left[\sbrac{0}, \, \sbrac{1/5}\right]\\
\gamma_2 & = & \left[\sbrac{0}, \, \sbrac{1/(2n+1)}\right]\\
\gamma_3 & = & \left[\sbrac{1}, \, \sbrac{1/2}\right]
\end{eqnarray*}
Note that $\gamma$ is not an admissible path system, but by adding vertical edges to the ends of the paths it extends to an admissible system.

\begin{lem} \label{lem:Kn_has_it}
For the knot $K_n$, the path system $\gamma$ extends to an admissible system that only carries incompressible surfaces with boundary slope $4(n+4)$.  Also, every path system that carries this slope is a vertical extension of $\gamma$. 
\end{lem}

\begin{proof} 
The cycle of final $r$-values for the path system $\gamma$ is $(1, \, -4, \, -2n, \, 1)$.  Furthermore, the final slopes of the $\gamma_i$ have positive sum, $\gamma$ satisfies (E2), and each $\gamma_i$ ends on the left edge of $\SS$.  Therefore, by Proposition 2.9 of \cite{HO}, the $\gamma_i$'s can be extended by vertical edges to form a system that carries an incompressible surface.  As in the proof of Claim 2 above, we can calculate the twist number of any such path system by calculating the twist number of the vertical extension $\gamma'$ given by adding the vertical edge $[\sbrac{0}, \sbrac{1}]$ to both $\gamma_0$ and $\gamma_3$.

We calculate the twist number as above; one has:
\begin{eqnarray*}
t(\gamma'_0) & = & 2(1-2) \ = \ -2, \\
t(\gamma'_1) & = & 2(1-0) \ = \ 2, \\
t(\gamma'_2) & = & 2(1-0) \ = \ 2, \\
t(\gamma'_3) & = & 2(1-1) \ = \ 0, 
\end{eqnarray*}
giving
$$t(\gamma') \ = \ \sum t(\gamma_i) \ = \ -2+2+2 \ = \ 2.$$
Therefore the boundary slope of any surface carried by a vertical extension of $\gamma$ is $t(\gamma')-t(s) = 2+14+4n=4(n+4)$. 

The uniqueness of $\gamma$ follows exactly as in Lemma \ref{lem:Ktau_doesn't}. 
\end{proof}

Combining Lemma \ref{lem:Ktau_doesn't} and Lemma \ref{lem:Kn_has_it} yields the following result: 

\begin{pro}
For every $n>1$, $4(n+4)$ is a boundary slope of $K_n$ but not of $K^\tau_n$. Furthermore, every path system that carries this slope is a vertical extension of $\gamma$. 
\end{pro}

\section{The slope $\mathbf{4(n+4)}$ of $\mathbf{K_n}$ is not strongly detected} 
\label{sec:strong}

The previous section has established that $4(n+4)$ is not a boundary slope of $K_n^\tau$. 
As stated in the introduction, if $4(n+4)$ is strongly detected by a curve in $X(K_n)$, then the restriction of every character to the Conway sphere $S_4$ must be reducible. This observation greatly simplifies the problem of showing that this boundary slope is not strongly detected.

\subsection{Some lemmata} 

 \begin{figure} [ht]
\begin{center}
\includegraphics[width=4.25in]{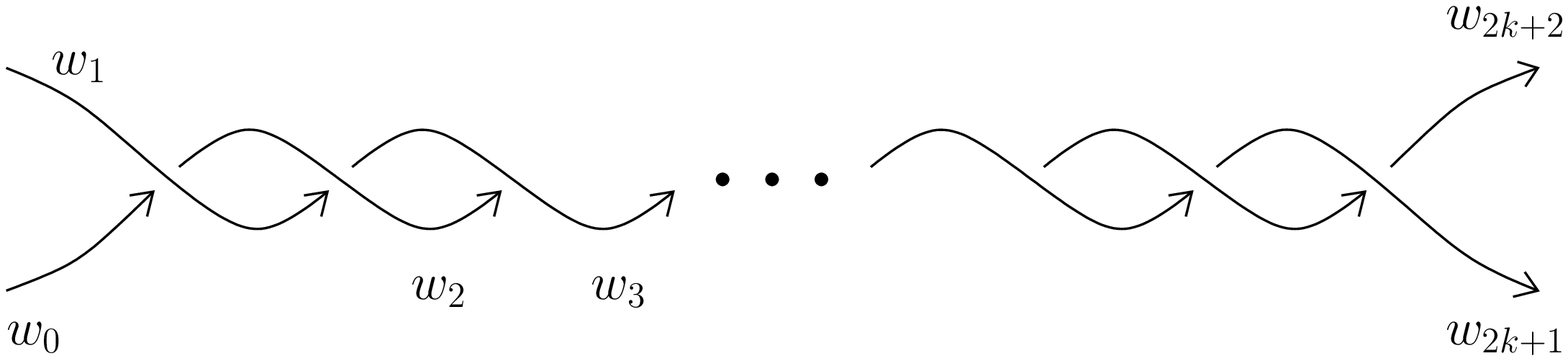}
\end{center}
\caption{The diagram for Lemma \ref{lem:twist}}  \label{fig:twist}
\end{figure}

\begin{lem} \label{lem:twist}
Assume that a twist region has a diagram as in Figure \ref{fig:twist} with Wirtinger generators $\{w_0, w_1, \ldots , w_{2k+2} \}$ as shown.  Then for every $k \geq 0$ we have
$$w_{2k+2} \ = \ (w_0 w_1)^{-k} w_1^{-1} (w_0 w_1)^{k+1}$$
and
$$ w_{2k+1} \ = \ (w_0 w_1)^{-k} w_1 (w_0 w_1)^{k}. $$
\end{lem}

\begin{proof}
For $k=0$, the relations are the Wirtinger relation and the trivial relation respectively. The conclusion follows by induction.
\end{proof}

\begin{lem} \label{lem:abelianize}
Assume that $x$ and $y$ are elements of a subgroup $G < \pie(M)$ where $x$ and $y$ have identical images in the abelianization of $G$.  If $\rho \in R(M)$ such that $\rho(G)$ is a group of upper triangular matrices, then $\rho(x)$ and $\rho(y)$ are identical along their diagonals.
\end{lem}

\begin{proof} Let $\Delta < \slc$ be the subgroup of upper-triangular matrices and $D < \slc$ be the abelian subgroup of diagonal matrices.  Then we have an epimorphism $\delta \co \Delta \longrightarrow D$ given by
$$\delta \left( \begin{array}{cc} \alpha & \beta \\ 0 & 1/\alpha \end{array} \right) \ = \ \left( \begin{array}{cc} \alpha & 0 \\ 0 & 1/\alpha \end{array} \right).$$
Since $x$ and $y$ have the same image in $G^{\mathrm{AB}}$ we have $\delta \rho(x)=\delta \rho(y).$
\end{proof}

\begin{lem} \label{lem:poly}
If $A=\left( \begin{smallmatrix} \alpha & 0 \\ 0 & 1/\alpha \end{smallmatrix} \right)$ and $B=\left( \begin{smallmatrix} \alpha & 1 \\ 0 & 1/\alpha \end{smallmatrix} \right)$ then for every $n \in \Z^+$
$$(AB)^n \ = \ \begin{pmatrix} \alpha^{2n} & p_n(\alpha) \\ 0 & \alpha^{-2n} \end{pmatrix}$$
where $p_n(x) \in \C(x)$ and $p_n(x) \neq \frac{1}{(1-x^2)x^{2n-1}}$.
\end{lem}

\begin{proof} 
First we establish an inductive formula for $p_n(x).$

\bigskip
{\bf Claim 1} $\qquad p_1(x) = x$ and $p_n(x) = x^{3-2n} + x^2 \, p_{n-1}(x)$.

For $n=1$ we calculate
$$AB \ = \ \begin{pmatrix} \alpha^2 & \alpha \\ 0 & \alpha^{-2} \end{pmatrix}$$ 
to see that $p_1(x)=x$ as claimed.  Now assume that 
$$(AB)^{n-1} \ = \ \begin{pmatrix} \alpha^{2n-2} & p_{n-1}(\alpha) \\ 0 & \alpha^{2-2n} \end{pmatrix}.$$
Multiplying this by $AB$ we see that 
$$(AB)^n \ = \ \begin{pmatrix} \alpha^{2n} & \alpha^2 \, p_{n-1}(\alpha) + \alpha^{3-2n} \\ 0 & \alpha^{-2n} \end{pmatrix}.$$
Hence, claim 1 is true.

\bigskip
{\bf Claim 2} $\qquad x^{2n-3} \, p_n(x) \in \C[x].$

Again, we establish this by induction.  For $n=1$ 
$$\inv{x} \, p_1(x) \ = \ \inv{x} \cdot x \ = \ 1 \ \in \C[x].$$
Now assume that $x^{2n-5} \, p_{n-1} (x) \in \C[x].$  Then, using claim 1, we have
\begin{eqnarray*}
 x^{2n-3} \, p_n(x) & = & x^{2n-3} \big(x^{3-2n} + x^2 \, p_{n-1}(x)\big) \\
 & = & 1 + x^{2n-1} \, p_{n-1}.
 \end{eqnarray*}
 This is in $\C[x]$ by the inductive assumption.
 
 Finally, we see that $p_n(x) \neq  \frac{1}{(1-x^2)x^{2n-1}}$ because 
 $$x^{2n-3} \cdot \frac{1}{(1-x^2)x^{2n-1}} \ = \ \frac{1}{(1-x^2)x^2} \ \notin \ \C[x]. $$
 \end{proof}
 
A Wirtinger presentation for $\pie(M)$ can be obtained from Figure \ref{fig:Wirt} with generating set $\{a,b,c,d,e,f\}$ given by the labels in the figure. We single out the following relations which will be used repeatedly in the arguments that follow:
\begin{itemize}
\item[(R1)] $ae=eb,$ and
\item[(R2)] $d=a\inv{b}c.$
\end{itemize}
Also, by applying Lemma \ref{lem:twist} to $T_0$ and $T_1$, we get
\begin{itemize}
\item[(R3)] $\inv{b} = fa\inv{f}\inv{a}\inv{f} \ ( \Leftrightarrow \, faf=bfa),$ and
\item[(R4)] $c=(fd)^{-2}\inv{d}(fd)^3 \ (\Leftrightarrow \, d(fd)^2c=(fd)^3).$
\end{itemize}

Before moving on to the calculations, we will prove one more elementary lemma.

\begin{lem} \label{lem:trace_zero}
If $\rho \in R(M)$ with $\rho(a) = \rho(\inv{b})$ and $\chi_\rho(a) \neq \pm 2$ then $\chi_\rho(a)=0.$
\end{lem}

\begin{proof}
Since $\chi_\rho(a) \neq \pm 2$ we may conjugate $\rho$ to assume that $\rho(a)$ is diagonal.  Let
$$\rho(a) \ = \ \begin{pmatrix} \alpha & 0 \\ 0 & \alpha \end{pmatrix} \AND \rho(e) \ = \ \begin{pmatrix} w & x \\ y & z \end{pmatrix}.$$
The lemma follows from the relation (R1), $ae=eb$.  We have
$$\rho(ae) \ = \ \begin{pmatrix}  \alpha w & \alpha x \\ y/\alpha & z/\alpha \end{pmatrix}$$
and
$$ \rho(eb) \ = \ \rho(e\inv{a}) \ = \ \begin{pmatrix}  w/\alpha  & \alpha x \\ y/\alpha & \alpha z \end{pmatrix}.$$
Then $\alpha w = w/\alpha.$  This implies that either $\alpha^2-1=0$ or $w=0$.  The first possibility is ruled out by assumption, hence $w=0$.  Similarly we have $z=0$.  Therefore $\chi_\rho(a) = \chi_\rho(e) = 0.$
\end{proof}

\subsection{Calculations}

The Conway sphere $S_4$ separates $M$ into two submanifolds $M_1$ and $M_2$.  Let $M_1$ be the piece that contains the tangles $T_0$ and $T_1$.  Choose a basepoint in $S_4$ for $\pie(M)$ and write $\Gamma_i=\pie(M_i)$ and $H=\pie(S_4)$. 

\begin{pro}
There exists a finite set $\Lambda \subset \C$ such that if $\rho$ is an irreducible representation in $R(M)$ where $\rho|_H$ is reducible, then $\chi_\rho(a) \in \Lambda$.  In particular, the boundary slope $4(n+4)$ is not strongly detected by $X(M)$.
\end{pro}

\begin{proof}
We will calculate all conjugacy classes of irreducible representations $\rho$ that are reducible on $H$ and show that the possible values for $\chi_\rho(a)$ lie in a finite set.  We start with $\Lambda = \{0, \pm 2 \}$ and add finitely many values to this set as we proceed through several cases.  

We may conjugate $\rho$ to assume that its restriction to $H$ is upper triangular.  That is, the matrices $\rho(a), \, \rho(b), \, \rho(c),$ and $\rho(d)$ are all upper triangular.  Furthermore, since we are looking for new values to add to $\Lambda$ we may assume that $\tr \big(\rho(a)\big) \notin \{0, \pm2\}$.

To simplify notation, given $g \in \pie(M)$, $\rho(g)$ will be denoted by $g$ for the remainder of this section.  The proof breaks into several cases and subcases.

\bigskip
{\bf Case 1}  $\quad \rho|_H$ is abelian.

\smallskip
Since $\tr(a) \neq \pm 2,$ we can conjugate $\rho$ to assume that $a$ is diagonal.  Since $\rho(H)$ is abelian we know that $b,c,d \in \{ a^{\pm1} \}$.  By Lemma \ref{lem:trace_zero}, we know that $b=a$.  Now using the relation (R1), we have $ae=ea$.  Hence $e$ is also diagonal.  Since $\Gamma_2$ is generated by $\{a,b,c,d,e \}$ we conclude that $\rho(\Gamma_2)$ is diagonal.  The representation $\rho$ is assumed to be irreducible, so $\rho|_{\Gamma_1}$ must be irreducible.

The relation (R2) implies that $d=c$ and so (R3) becomes $faf=afa$ and (R4) becomes $cfcfc=fcfcf$.

We have already established that $c \in \{a^{\pm 1} \}$.  
\begin{itemize}
\item Subcase A $\quad c=a.$

The relation (R4) becomes $afafa=fafaf$.  Using (R3) we have
\begin{eqnarray*}
afafa & = & fa(faf) \\
& = & fa(afa) \\
af & = & fa.
\end{eqnarray*}
Then $f$ is diagonal which contradicts that $\rho|_{\Gamma_1}$ is irreducible.
\item Subcase B $\quad c=\inv{a}.$

Since $\rho|_{\Gamma_1}$ is irreducible we can reconjugate the representation to assume that
$$a \ = \ \begin{pmatrix}  \alpha & 1 \\ 0 & 1/\alpha \end{pmatrix} \AND f \ = \ \begin{pmatrix}  \beta & 0 \\ r & 1/\beta \end{pmatrix}$$
where $\beta \in \{ \alpha^{\pm 1} \}$ and $r \neq 0$.

If $\beta = \alpha$ then (R3) implies that $1=\alpha^2 +r +\frac{1}{\alpha^2}$ and (R4) implies that $r^2-3r+1=0.$  There are at most eight simultaneous solutions to these two equations.  We add the corresponding values of $\tr(a)$ to $\Lambda$ and note that $\Lambda$ remains finite.

Otherwise, we have $\beta=1/\alpha$.  This time (R3) forces $r=-1$.  Together with (R4), this means that $\alpha$ must satisfy the equation $\alpha^8+\alpha^6+\alpha^4+\alpha^2+1=0$.  Again, we need only add a finite number of values to $\Lambda$ to account for these representations.
\end{itemize}    

We may now assume that $\rho(H)$ is not abelian.

\bigskip
{\bf Case 2} $\quad \rho|_{\Gamma_2}$ is reducible.

\smallskip
In this case, we may conjugate $\rho$ to assume that $a$ is diagonal
and $b, \, c, \, d,$ and $e$ are all upper triangular.
\begin{itemize}
\item Subcase A $\quad [a,b]=I.$

Then $b \in \{a^{\pm 1} \}.$  Lemma \ref{lem:abelianize} shows that we must have $b=a.$  Then (R1) gives $[a,e]=I$ and $e$ must also be diagonal.  Also, (R2) implies $d=c$.  Since $\rho(H)$ is non-abelian, $c$ cannot be diagonal.

Let 
$$e \ = \ \begin{pmatrix} \alpha & 0 \\ 0 & 1/\alpha \end{pmatrix}.$$
and by conjugating by diagonal matrices if necessary we may assume that
$$c \ = \ \begin{pmatrix} \beta & 1 \\ 0 & 1/\beta \end{pmatrix}$$
where $\beta \in \{ \alpha^{\pm 1} \}.$   Note that in $\Gamma_2^{\mathrm{AB}}$ we have $c=e$.  Since $\rho|_{\Gamma_2}$ is reducible we may apply Lemma \ref{lem:abelianize} to see that $c$ and $e$ are identical on their diagonals.  Hence $\beta= \alpha$.

Using Lemma \ref{lem:twist}, we have $(ec)^n e= c (ec)^n$.  Also, by Lemma \ref{lem:poly} we have$$ (ec)^n \ = \ \begin{pmatrix} \alpha^{2n} & p_n(\alpha) \\ 0 & \alpha^{-2n} \end{pmatrix}$$
where $p_n(x) \in \C(x)$ and $p_n(x) \neq \frac{1}{(1-x^2)x^{2n-1}}$.

Equating the upper right entries of the matrices $(ec)^n e$ and $c(ec)^n$ we get 
$$ \alpha \, p_n(\alpha) + \alpha^{-2n} \ = \ \frac{1}{\alpha} \, p_n(\alpha).$$  
We have assumed that $\alpha \notin \{0, \pm 1 \}$ so we have
$$ p_n(\alpha) \ = \ \frac{1}{(1-\alpha^2)\alpha^{2n-1}}.$$
This equation has finitely many solutions because $p_n(x) \neq \frac{1}{(1-x^2)x^{2n-1}}$.  We add the corresponding values $\alpha + \inv{\alpha}$ to $\Lambda$ and move on.

\item Subcase B $\quad [a,b] \neq I.$

Recall that we have conjugated $\rho$ so that $a$ is diagonal and $\rho(\Gamma_2)$ is upper triangular.  Since $[a,b] \neq I$ we know that the upper right entry of $b$ is non-zero, so we can conjugate by diagonal matrices to assume this entry is $1$.  As before, Lemma \ref{lem:abelianize} implies that $a$ and $b$ are identical along their diagonals.  We have
$$a \ = \ \begin{pmatrix} \alpha & 0 \\ 0 & 1/\alpha \end{pmatrix} \AND b \ = \ \begin{pmatrix} \alpha & 1 \\ 0 & 1/\alpha \end{pmatrix}.$$
We have assumed that $\rho$ is irreducible.  This implies that the lower left entry of $f$ is non-zero.  We also know that $\tr(f) = \tr(a) = \alpha + 1/\alpha$.  Putting this together, we know that $f$ is of the form
$$ f \ = \ \begin{pmatrix} x & \frac{1}{\alpha y} \, (-\alpha x^2  + \alpha^2 x + x - \alpha) \\ y & \frac{1}{\alpha} \, (-\alpha x + \alpha^2 +1) \end{pmatrix}$$
where $y \neq 0.$  

If $\rho$ is a representation then we must have $faf=bfa$ by (R3).  Looking at the lower left entries of these matrices we get the equation
$$\alpha x y + \frac{y}{\alpha^2} \, (-\alpha x + \alpha^2 +1).$$
Then because $y\neq 0$ we must have
$$x \ = \ \frac{-1}{\alpha(\alpha^2-1)}.$$
Making this substitution and equating the upper left entries of $faf$ and $bfa$, we arrive at the conclusion that $y\alpha =0$, which is not possible.  Therefore, there are no such representations.

 \bigskip
{\bf Case 3} $\rho|_{\Gamma_2}$ is irreducible.

\smallskip
In this case, when we conjugate so that $a$ is diagonal and $b, c,$ and $d$ are upper triangular, we know that the lower left entry of $e$ must be non-zero.

First, note that $[a,b] \neq I$.  This follows from (R1) and Lemma \ref{lem:trace_zero}.  Since $a$ is diagonal and not parabolic, if $a$ and $b$ commute then $b=\inv{a}$ or $b=a$.  The lemma shows that the first possibility doesn't occur.  On the other hand, if $a=b$, then (R1) implies that $a$ and $e$ commute which is untrue because the lower left entry of $e$ is non-zero.

We may now assume that 
$$a \ = \ \begin{pmatrix} \alpha & 0 \\ 0 & 1/\alpha \end{pmatrix} \AND b \ = \ \begin{pmatrix} \beta & 1 \\ 0 & 1/\beta \end{pmatrix}$$
where $\beta \in \{ \alpha^{\pm1}\}.$  Also, as argued for $f$ in Case 2: Subcase B, we have
$$ e \ = \ \begin{pmatrix} x & \frac{1}{\alpha y} \, (-\alpha x^2  + \alpha^2 x + x - \alpha) \\ y & \frac{1}{\alpha} \, (-\alpha x + \alpha^2 +1) \end{pmatrix}$$
where $y \neq 0.$   

We claim that $\beta \neq \alpha$.  This follows from (R1).  Equating the lower right entries of $ae$ and $eb$ we get
$$\frac{1}{\alpha^2} \, (-\alpha x + \alpha^2 +1) \ = \ y + \frac{1}{\alpha^2} \, (-\alpha x + \alpha^2 +1)$$
which implies that $y=0$, a contradiction.  Thus we must have 
$$b \ = \ \begin{pmatrix} 1/\alpha & 1 \\ 0 & \alpha \end{pmatrix}.$$

Looking at the upper right entries of $ae$ and $eb$, we see that $x=0$.  Substituting for $x$ and checking $ae=eb$ in the lower right entries, we conclude that $y=\alpha^{-2} (1-\alpha^4)$.  This implies that
$$e \ = \ \begin{pmatrix} 0 & \frac{\alpha^2}{\alpha^4 -1} \\ \frac{1-\alpha^4}{\alpha^2} & \frac{\alpha^2+1}{\alpha} \end{pmatrix}.$$ 

We know that $c$ is upper triangular and $\tr(c)=\tr(a)$, thus 
$$c \ = \ \begin{pmatrix} \beta & r \\ 0 & 1/\beta \end{pmatrix}$$
where $\beta \in \{\alpha^{\pm 1} \}$ and $r \in \C$.  Formally assign the matrix $\Theta \in \slc$ to represent the product $(ec)^n$. 

Applying Lemma \ref{lem:twist} to the tangle $T_2,$ we get the relation $w_{2n-1} = (ec)^{-n}c(ec)^n$ and from $T_3$ we have that $w_{2n-1} = ae\inv{a}$.  Putting these together, we conclude that the matrix $\Theta$ should satisfy the relation 
$$\inv{\Theta}c\,\Theta=ae\inv{a}. \quad \mathrm{(R5)}$$
The matrix $\Theta$ should commute with $ec$ so we also have
$$\Theta \, ec = ec\, \Theta. \quad \mathrm{(R6)}$$

Using (R5) and $\det(\Theta)=1$, we can express $\Theta$ in terms of $\alpha, \, \beta, r,$ and a new independent variable $z$.  By comparing the matrices in (R6), we can express $r$ as a function in $\alpha, \, \beta,$ and $z$.  If $\beta = \alpha$ then the upper right entries of the matrices from (R5) gives the equality $\alpha^6 - \alpha^4 - \alpha^2 + 1 = 0$.  Thus we may assume that $\beta=1/\alpha$ at the 
price of adding  $\{\, \alpha + \inv{\alpha} \, | \, \alpha^6 - \alpha^4 - \alpha^2 + 1 = 0 \, \}$ to $\Lambda$.  Both $c$ and $\Theta$ are now expressed as functions of $\alpha$ and $z$:
$$c \ = \ \begin{pmatrix} 1/\alpha & r \\ 0 & \alpha \end{pmatrix} \AND \Theta \ = \ \begin{pmatrix} \frac{\alpha^4-z^2}{\alpha^4 z} & \frac{\alpha z}{\alpha^4-1} \\ \frac{z(1-\alpha^4)}{\alpha^5} & z \end{pmatrix}$$
where $r=\frac{\alpha^4 + \alpha^2z^2 - z^2}{z^2(\alpha^4 - 1)}$.

Write
$$f \ = \ \begin{pmatrix} p     & q \\ s & t \end{pmatrix}$$
and consider the consequences of (R3).  Equating the lower right entries of $faf$ and $bfa$,  we see that $q(\alpha^3 p + \alpha t - 1) - \alpha t = 0.$  If $q=0$ this equation implies that $t=0$. Then $\det(f) = 0,$ a contradiction.  Hence $q \neq 0$ and we solve $\det(f) = 1$ for $s$ to get $s = \frac{pt-1}{q}$.  Substituting this for $s$ and looking at the lower right entries again, we get an expression for $p$ in terms of $\alpha$ and $t$.  Compare the lower left entries to see that $(\alpha - t)(\alpha^2 t - t - \alpha)=0$.  There are two subcases.    

\begin{itemize}
\item[$\bullet$] Subcase A $\quad t = \alpha.$  \hspace{.25in} (This is the case where $\rho|_{\Gamma_1}$ is reducible.)

At this point we have $f$ in terms of $\alpha$ and $q$.  Using (R3) once again, we are able to solve for $q$ to get
$$ f \ = \ \begin{pmatrix} 1/\alpha & \frac{\alpha^2}{2 \alpha^2 -1} \\ 0 & \alpha \end{pmatrix}.$$
Checking (R4) gives the polynomial equation
\begin{equation*} 
(6 \alpha^8 -17 \alpha^6 +13 \alpha^4 +2 \alpha^2 - 4)z^2 - 6 \alpha^8 + 7 \alpha^6 - 2 \alpha^4 \ = \ 0.
\end{equation*}
By adding 
$\{ \alpha + \inv{\alpha} \, | \, 6 \alpha^8 -17 \alpha^6 +13 \alpha^4 +2 \alpha^2 - 4 \, = \, 0 \}$ to the set $\Lambda,$ we can solve for $z^2$ to get
\begin{equation*} 
z^2 \ = \ \frac{\alpha^4\,(6 \alpha^4 - 7 \alpha^2 +2)}{6 \alpha^8 - 17 \alpha^6 + 13 \alpha^4 + 2 \alpha^2 - 4}.
\end{equation*}
We would like to use this equality to replace $z$ in our only remaining relation, $\Theta = (ec)^n$.
Although our expression for $\Theta$ in terms of $\alpha$ and $z$ does contain odd powers of $z$, when we look at the resulting expression for $\Theta^2$ in terms of $\alpha$ and $z$ we see that it contains only even powers of $z$.  We replace our expression for $z^2$ into $\tr(ec)$ and $\tr(\Theta^2)$ to write these quantities as rational functions in $\alpha$ only,
$$\tr(ec) \ = \ \frac{(\alpha^2+1)(10\alpha^4-15\alpha^2+6)}{\alpha^4\,(2\alpha^2-1)(3\alpha^2-2)}$$
and
$$\tr(\Theta^2) \ = \ \frac{\scriptstyle (\alpha^2+1)(72\alpha^{12}-288\alpha^{10}+446\alpha^8-278\alpha^6-25\alpha^4+108\alpha^2-36)}{\scriptstyle \alpha^4\,(2\alpha^2-1)(3\alpha^2-2)(6\alpha^6-11\alpha^4+2\alpha^2+4)}.$$

In light of this, instead of using $\Theta=(ec)^n$, we focus on the relation 
\begin{equation} \label{eq:reducible}
\tr(\Theta^2) = \tr\big((ec)^{2n}\big).
\end{equation}
Using standard trace formulas, $\tr(\Theta^2) = P\big(\tr(ec)\big)$ where $P(x) \in \Z[x]$ is a monic polynomial of degree $2n$.  We wish to rule out the possibility that there are infinitely many solutions for $\alpha$ that will satisfy this equation.  This will happen precisely when the rational functions $\tr(\Theta^2)$ and $P\big(\tr(ec)\big)$ are identically equal.  However, $\tr(\Theta^2)$ has a pole at the roots of the polynomial $6x^6-11x^4+2x^2+4$ while $P\big(\tr(ec)\big)$ takes a finite value at each of these roots, since the poles of $P\big(\tr(ec)\big)$ are the same as the poles of $\tr(ec)$.  Therefore, these rational functions are not identically equal and there are only finitely many values for $\alpha$ that satisfy equation (\ref{eq:reducible}).  We add the corresponding traces to $\Lambda$ and move to the second subcase.  

\bigskip
\item[$\bullet$] Subcase B $\quad t = \frac{\alpha}{\alpha^2-1}.$ \hspace{.25in} ($\rho|_{\Gamma_1}$ is irreducible.)

Just as in subcase A, we can use (R3) to get
$$f \ = \ \begin{pmatrix} \frac{\alpha^4-\alpha^2-1}{\alpha(\alpha^2-1)} & \frac{\alpha^2}{(\alpha^2+1)(\alpha^2-1)^2} \\ \frac{\alpha^4 -\alpha^2-2}{\alpha^2} & \frac{\alpha}{\alpha^2-1} \end{pmatrix}.$$
After adding $\pm\big(2^{1/2} + 2^{-1/2}\big)$ to $\Lambda$, the relation (R4) implies that 
\begin{multline} \label{eq:irred}
(-\alpha^{12} + 8 \alpha^{10} -22 \alpha^8 +23 \alpha^6-3 \alpha^4 - 9 \alpha^2 +4)z^4 + \\
\quad + (2\alpha^{12}-12\alpha^{10}+22\alpha^8-11\alpha^6-3\alpha^4+2\alpha^2)z^2- \\ \quad -\alpha^{12}+4\alpha^{10}-4\alpha^8 \ = \ 0.
\end{multline}
Let $\cc \subset \C^2$ be the zero set for this polynomial.  Setting the variables 
\begin{equation} \label{eq:T}
T \ = \ \tr(\Theta) \ = \ \frac{\alpha^4z^2+\alpha^4-z^2}{\alpha^4z}
\end{equation}
and
\begin{equation} \label{eq:E}
E \ = \ \tr(ec)\ = \ \frac{\alpha^4z^2-\alpha^4+z^2}{\alpha^2z^2}
\end{equation}
we have a corresponding projection 
$$\phi \co \cc \longrightarrow \big\{ \, (T,E) \in \C^2 \, \big\}.$$
We take resultants to get a polynomial $Q(E,T) \in \Z[E,T]$, where $Q=0$ on the image of $\phi.$  (This polynomial is quite large, so we do not include it here.)  As before, we can express $T$ as a monic degree $n$ polynomial $r(E) \in \Z[E]$.  The one variable polynomial equation $Q\big(E,r(E)\big) = 0$ holds on $\mathrm{Im}(\phi)$.  As long as this polynomial is not identically zero, the image of $\phi$ is $0$-dimensional.   We show that this is true by showing that its degree is bigger than zero.

Writing $Q$ as a polynomial in $\Z[E][T]$, we get
$$ Q(E,T) \ = \ \sum^5_{i=0} \ p_i(E)\cdot T^{2i}$$
where 
\begin{eqnarray*}
\deg\big(p_0(E)\big) \ = \ 12 & \Longrightarrow & \deg\big(p_0(E) \cdot r(E)^0\big) \ = \ 12 \\ 
\deg\big(p_1(E)\big) \ = \ 11 & \Longrightarrow & \deg\big(p_1(E) \cdot r(E)^2\big) \ = \ 11+2n \\  
\deg\big(p_2(E)\big) \ = \ 10 & \Longrightarrow & \deg\big(p_2(E) \cdot r(E)^4\big) \ = \ 10+4n\\ 
\deg\big(p_3(E)\big) \ = \ 8 & \Longrightarrow & \deg\big(p_3(E) \cdot r(E)^6\big) \ = \ 8+6n \\ 
\deg\big(p_4(E)\big) \ = \ 7 & \Longrightarrow & \deg\big(p_4(E) \cdot r(E)^8\big) \ = \ 7+8n \\ 
\deg\big(p_5(E)\big) \ = \ 4 & \Longrightarrow & \deg\big(p_5(E) \cdot r(E)^{10}\big) \ = \ 4 + 10n. 
\end{eqnarray*}
Since $n>1$, the degree $4+10n$ term in $p_5(E)\cdot r(E)^{10}$ is the unique term of highest degree in $Q\big(E,r(E)\big)$.  Therefore,  
$$\deg\big(Q\big(E,r(E)\big) \ = \ 4+10n \  > 0.$$

We have now established that $\mathrm{Im}(\phi)$ is a finite set of points.  The proof will be complete if we show that the fibers of $\phi$ are finite.  Using equations (\ref{eq:T}) and (\ref{eq:E}),  
$$z \ = \ \frac{\alpha^2 T}{2\alpha^2-E}.$$
This, together with equation (\ref{eq:E}), implies that
$$0 \ = \ (T^2-4)\alpha^4+(4E-T^2E)\alpha^2-E^2+T^2.$$
Therefore, the fibers are finite over every $(E,T)$ unless $T^2=4$ and $E^2=4$.  It is easy to see that these points are not in the image of $\phi.$ 
\end{itemize}
\end{itemize}
This completes the proof that $4(n+4)$ is not strongly detected.
\end{proof}

\bibliographystyle{plain}
\bibliography{Biblio}


\address{Department of Mathematics, The University of Texas at Austin, Austin, TX 78712, USA}
\email{chesebro@math.utexas.edu}

\address{D\'epartment de math\'ematiques,
  Universit\'e du Qu\'ebec \`a Montr\'eal,
  Case postale 8888, Succursale Centre-Ville,
  Montr\'eal (Qu\'ebec) Canada H3C 3P8} 
\email{tillmann@math.uqam.ca} 
\Addresses
\end{document}